\documentclass{amsart}

\usepackage{epsfig}
\usepackage{latexsym}
\theoremstyle{plain}

\newtheorem{thm}{Theorem}

\usepackage[all]{xy}

\theoremstyle{remark}

\usepackage{amssymb}
\usepackage{amsthm}
\usepackage{amsmath}
\usepackage{amsxtra}

\DeclareMathOperator{\Bez}{Bez}

\newcommand{\bfP}{\mathbb P}
\newcommand{\bfR}{\mathbb R}

\newcommand{\Q}{\mathbb Q}

\newcommand{\Z}{\mathbb Z}
\newcommand{\N}{\mathbb N}

\newcommand{\LL}{{\mathcal L}}

{.tfm}
{.tfm}

\begin{document}

\title{
Integral ratios of factorials and algebraic hypergeometric functions
}

\author{ Fernando Rodriguez-Villegas}
\address{Department of Mathematics 
University of Texas at Austin, TX 78712}

\email {villegas\@math.utexas.edu}

\maketitle

Chebychev in his work on the distribution of primes numbers used the
following fact
\[
u_n:=\frac{(30n)!n!}{(15n)!(10n)!(6n)!} \in \Z, \qquad n=0,1,2,\ldots
\] This is not immediately obvious (for example, this ratio of
factorials is not a product of multinomial coefficients) but it is not
hard to prove. The only proof I know proceeds by checking that the
valuations $v_p(u_n)$ are non-negative for every prime $p$; an
interpretation of $u_n$ as counting natural objects or being
dimensions of natural vector spaces is far from clear.

As it turns out, the generating function
\[
u:=\sum_{\nu\geq 1}u_n\lambda^n
\]
is algebraic over $\Q(\lambda)$; i.e. there is a polynomial $F\in
\Z[x,y]$ such that
\[
F(\lambda,u(\lambda))=0.
\]
However, we are not likely to see this polynomial explicitly any time
soon as its degree is $483,840$ (!)

What is the connection between $u_n$ being an integer for all $n$ and
$u$ being algebraic? Consider the more general situation
\[
u_n:=\prod_{\nu\geq 1} (\nu n)!^{\gamma_\nu},
\]
where the sequence $\gamma=(\gamma_\nu)$ for $\nu \in \N$ consists of
integers which are zero except for finitely many.

We assume throughout that $\gamma$ is {\it regular}, i.e.,
\[
\sum_{\nu \geq 1} \nu \gamma_\nu =0,
\]
 which, by Stirling's formula, is equivalent to the generating
series $u:=\sum_{\nu\geq 1}u_n\lambda^n$ having finite non-zero radius
of convergence. We define the {\it dimension} of $\gamma$ to be
\[
d:=-\sum_{\nu\geq 1} \gamma_\nu.
\]
 To abbreviate, we will say that $\gamma$ is {\it integral} if
$u_n\in \Z$ for every $n=0,1,2,\ldots$.

We can now state the main theorem of the talk.

\begin{thm}
Let $\gamma\neq 0$ be regular; then $u$ is algebraic if and only if
$\gamma$ is integral and $d=1$.
\end{thm}

One direction is fairly straightforward. If $u$ is algebraic, by a
theorem of Eisenstein, there exists an $N\in \N$ such that $N^nu_n \in
\N$ for all $n\in \N$. It is not hard to see that in our case if such
an $N$ exists then it must equal $1$. To see that $d=1$ we need to
introduce the {\it monodromy representation}. 

The power series $u$ satisfies a linear differential equation $Lu=0$.
After possibly scaling $\lambda$ this equation has singularities only
at $0,1$ and $\infty$. Indeed, $u$ is a hypergeometric
series. Moreover, these singularities are regular singularities
precisely because we assumed $\gamma$ to be regular.

If we let $V$ be the space of local solutions to $Lu=0$ at some base
point not $0,1$ or $\infty$ then analytic continuation gives a
representation
\[
\rho: \pi_1(\bfP^1\setminus \{0,1,\infty\}) \longrightarrow GL(V).
\] We let the {\it monodromy group} $\Gamma$ be the image of $\rho$
and let $B,A,\sigma$ be the monodromies around $0,\infty,1$,
respectively, with orientations chosen so that $A=B\sigma$. The main
use of the monodromy group for us is the fact that $u$ is algebraic if
and only if $\Gamma$ is finite.

As it happens the multiplicity of the eigenvalue $1$ for $B$ is $d$
and it is also true that the corresponding Jordan block of $B$ is of size
$d$. Hence, $\Gamma$ is not finite if $d>1$.

To prove the converse we appeal to the work of Beukers and Heckman
\cite{BH} who extended Schwartz work and described all algebraic
hypergeometric functions. Let $p$ and $q$ be the characteristic
polynomials of $A$ and $B$ respectively. In our situation $p$ and $q$
are relatively prime polinomials in $\Z[x]$ (which are products of
cyclotomic polynomials). Their work tells us that $\Gamma$ is finite
if and only if the roots of $p$ and $q$ interlace in the unit circle.

The key step in the proof of this beautiful fact is to determine when
$\Gamma$ fixes a non-trivial positive definite Hermitian form $H$ on
$V$ (which guarantees that $\Gamma$ is compact). I explained in my
talk how $H$ can be defined using a variant of a construction going
back to Bezout. Consider the two variable polynomial
\[
\frac{p(x)q(y)-p(y)q(x)}{x-y} =\sum_{i,j}B_{i,j}x^iy^k
\]
and define the  {\it Bezoutian} of $p$ and $q$ as 
\[
\Bez(p,q)=(B_{i,j}).
\]
 We need two facts about this matrix. First, the determinant of
$\Bez(p,q)$ equals the resultant of $p$ and $q$ (in passing I should
mention that this is a useful fact computationally since the matrix is
of smaller size than the usual Sylvester matrix). Second, note that
$\Bez(p,q)$ is symmetric. Hence it carries more information than just
its determinant as it defines a quadratic form $H$. It is a classical
fact (due to Hermite and Hurwitz) that the signature of $H$ has a
topological interpretation.

 Consider the continuous map $\bfP^1(\bfR)\rightarrow \bfP^1(\bfR)$
given by the rational function $p/q$. Since $\bfP^1(\bfR)$ is
topologically a circle we have $H^1(\bfP^1(\bfR),\Z) \simeq \Z$ and
the induced map $H^1(\bfP^1(\bfR),\Z) \rightarrow
H^1(\bfP^1(\bfR),\Z)$ is multiplication by some integer $s$, which is
none other than the signature of $H$. In particular, $H$ is definite
if and only if the roots of $p$ and $q$ interlace on $\bfR$.  A
twisted form of this construction and analogous signature result can
be applied to the hypergeometric situation; in this way we recover the
facts about the Hermitian form fixed by $\Gamma$ proved by Beukers and
Heckman.

Finally, to make the connection with the integrality of $\gamma$ we
define the {\it Landau function}
\[
\LL(x) :=-\sum_{\nu \geq 1} \gamma_\nu \{\nu x\}, \qquad
x \in \bfR  
\]
where $\{x\}$ denotes  fractional part. It
is simple to verify that
\[
v_p(u_n)=\sum_{k\geq 1} \LL\left(\frac n{p^k}\right).
\]
 Landau \cite{L} proved a nice criterion for integrality: $\gamma$ is
integral if and only if $ \LL(x)\geq 0$ for all $ x \in \bfR$. 

Write 
\[
p(t)=\prod_{j=1}^r(t-e^{2\pi i\alpha_j}), \qquad
q(t)=\prod_{j=1}^r(t-e^{2\pi i\beta_j}),
\]
where $r=\dim V$ and  $0\leq\alpha_1\leq \alpha_2\leq \cdots \leq
\alpha_r<1$ and $0\leq \beta_1  \leq \beta_2 \leq\cdots \leq
\beta_r<1$ are  rational. 

The function $\LL$ satisfies a number of simple properties: it is
locally constant (by regularity), periodic modulo $1$, right
continuous with discontinuity points exactly at $x\equiv \alpha_j
\bmod 1$ or $x\equiv \beta_j \bmod 1$ for some $j=1,\ldots, r$ and
takes only integer values.  More precisely,
\[
\LL(x)=\#\{j  \;|\;\alpha_j\leq x\}-\#\{j\;|\;0<\beta_j\leq x\}.
\]
Away from the discontinuity points of $\LL$ we have
\[
\LL(-x)=d-\LL(x).
\]
In particular, $\LL(x)\geq 0$ if and only if $\LL(x)\leq d$.

It is now easy to verify that if $d=1$ and $\LL(x)\geq 0$ then the
 roots of $p$ and $q$ must necessarily interlace on the unit circle
 finishing the proof. (Some further elaboration would also yield the
 other implication in the theorem independently of our previous
 argument.)

As a final note, let me mention that the examples in the theorem are a
case of the ADE phenomenon; up to the obvious scaling $n \mapsto dn$
for some $d\in \N$, they come in two infinite families $A$ and $D$,
which are easy to describe, and some sporadic ones ($10$ of type
$E_6$, $10$ of type $E_7$ and $30$ of type $E_8$).

\end{document}